\documentclass[conference,10pts]{IEEEtran}

\ifCLASSINFOpdf

\else

\fi

\usepackage{amssymb,amsmath}
\usepackage{amsmath}
\usepackage{amsfonts}
\usepackage{MnSymbol,enumerate}

\ifCLASSINFOpdf

\else

\fi

\usepackage{amssymb}
\usepackage[dvips,final]{epsfig}
\usepackage{amsfonts}
\usepackage{latexsym}

\usepackage{color}

\usepackage{setspace}

\newcommand{\oo}{\mbox{$\mathbb O$}}

\newcommand{\uu}{\mbox{$\mathbf u$}}
\newcommand{\x}{\mbox{$\mathbf x$}}
\newcommand{\y}{\mbox{$\mathbf y$}}
\newcommand{\vv}{\mbox{$\mathbf v$}}

\newcommand{\ttt}{\mbox{$\cal T$}}
\newcommand{\f}{\mbox{$\cal F$}}

\newcommand{\qa}{\mbox{\quad\mbox{and}\quad}}

\newcommand{\rt}{\mbox{$\mathbb R$}}

 \newcommand{\rank}{\mathrm{rank\;}}

\newcommand{\bxi}{\mbox{\boldmath $\xi$}}

 \def\diag{\mathop{{\rm diag}}\nolimits}

\newtheorem{theorem}{Theorem}

\newtheorem{example}{Example}

\newtheorem{definition}{Definition}
\newtheorem{remark}{Remark}

\begin{document}

\title{Data Compression: Multi-Term  Approach}

\author{Pablo Soto-Quiros\IEEEauthorrefmark{1}$^,$\IEEEauthorrefmark{2} and \IEEEauthorblockN{Anatoli Torokhti\IEEEauthorrefmark{1}}
\IEEEauthorblockA{\IEEEauthorrefmark{1}Centre for Industrial and Applied Mathematics, University of South Australia, Adelaide, SA 5095, Australia}
\IEEEauthorblockA{\IEEEauthorrefmark{2}Escuela de Matem\'aticas, Instituto Tecnol\'ogico de Costa Rica, Cartago 30101, Costa Rica\\
Emails: \{juan.soto-quiros, anatoli.torokhti\}@unisa.edu.au}}


%
\maketitle
\begin{abstract}

 In terms of signal samples, we propose  and justify a new rank reduced multi-term transform, abbreviated as MTT,  which, under certain conditions, may  provide the better associated accuracy than that of known optimal rank reduced  transforms. The basic idea is to construct the transform with more parameters to optimize than those in the known optimal transforms.
 This is realized by the extension of the known transform structures to the form that includes additional terms - the MTT has four matrices to minimize the cost. The MTT structure has also a special transformation that decreases the numerical load. As a result, the MTT performance is improved by the variation of the MTT components.

\end{abstract}

\section{Introduction}

Data compression is one of the central and widely studied problems in communication and signal processing. Assume, a reference signal  and  a noisy observed signal are represented by random vectors $\x\in  L^2(\Omega,\mathbb{R}^m)$ and $\y\in  L^2(\Omega,\mathbb{R}^n)$, respectively\footnote{Here, $\Omega=\{\omega\}$ is the set of outcomes, $\Sigma$ a $\sigma$-field of measurable subsets of $\Omega$, $\mu:\Sigma\rightarrow[0,1]$ an associated probability measure on $\Sigma$ with $\mu(\Omega)=1$ and $(\Omega,\Sigma,\mu)$ for a probability space.}. Then $\y$ compression (and denoising)  to a `shorter' (or compressed) vector $\uu\in  L^2(\Omega,\mathbb{R}^k)$ where $k \leq \min \{m, n\},$ and its consequent reconstruction (or decompression) can be represented as $F=DC$, where $C$ and $D$ are a `compressor' and `de-compressor', respectively.

An optimal determination of the  compressor  and de-compressor in the form of matrices  $C\in\rt^{k\times n}$  and  $D\in\rt^{m\times k}$, respectively, is given by the transform by Brillinger (BT) \cite{Brillinger2001},  the generalized Brillinger transform (GBT1) \cite{tor843} and the generic Karhunen-Lo\`{e}ve transform (GKLT) \cite{681430,torbook2007}. Under some  assumptions, the  BT, GKLT and GBT1 provide the best  associated  accuracy  among all linear transforms of the form $F=DC$,  for the same compression ratio $\displaystyle c=k/{\min \{m, n\}}$, where $k=1,2,\ldots, \min \{m, n\}$.
The second degree transforms developed in \cite{tor843,2001309} (and abbreviated here as the GBT2 and GKLT2, respectively), under the certain condition, improve the GKLT and GBT1 accuracy, for the same compression ratio. Nevertheless, it may happen that the accuracy is still not satisfactory.

 In the practical setting or in the training stage, signals are replaced with their samples. Here, we consider the case.
 Sample matrices of ${\x}$ and ${\y}$ are denote by $X\in\mathbb{R}^{m\times s}$ and $Y\in\mathbb{R}^{n\times s}$,   respectively, where $s$ is a number of samples.

 In terms of the signal samples, we propose  and justify a new rank reduced multi-term transform, abbreviated as MTT,  which, under certain conditions, may  provide the better associated accuracy than that of the  GBT1, GKLT, GKLT2 and GBT2.
 This is achieved by the extension of the known transform structures to the form that includes additional terms (see Section \ref{new}). As a result, the MTT accuracy is improved by the variation of the MTT components.  See Sections \ref{bdc1} and \ref{nma9}, and Examples \ref{ui289}-\ref{vnm82} below.

\subsection{Preliminaries}


Each matrix $T\in\rt^{m\times n}$ defines a bounded linear transformation $\ttt: L^2(\Omega,\mathbb{R}^n) \rightarrow L^2(\Omega,\mathbb{R}^m)$. We write $T$ rather then $\ttt$ since $[\ttt(\x)](\omega) = T[\x(\omega)]$, for each $\omega\in \Omega$.

 The covariance matrix formed from $\x$ and $\y$ is denoted by $E_{xy}.$ The pseudo-inverse  of matrix $M$ is denoted by $M^{\dag}$.

Let the SVD of matrix  $A\in \rt^{m\times n}$ be given by  ${A=U_A\Sigma_A V_A^T}$, where  ${U_A=[u_1 \;u_2\;\ldots u_m]\in \rt^{m\times m}}$, ${V_A=[v_1 \;v_2\;\ldots v_n]\in \rt^{n\times n}}$  are unitary matrices,  ${\Sigma_A=\diag(\sigma_1(A),\ldots,\sigma_{\min(m,n)}(A))\in\rt^{m\times n}}$ and $\sigma_1(A)$ $\ge \sigma_2(A)\ge\ldots\ge 0$. For $k<m$, $j<n$ and $\ell<\min(m,n)$, we denote $U_{A,k}:=[u_1 \;u_2\;\ldots u_k]$, $V_{A,j}=[v_1 \;v_2\;\ldots v_j]$ and $\Sigma_{A,\ell}=\diag(\sigma_1(A), \ldots,\sigma_\ell(A))$. Then
$\displaystyle [A]_k= \sum_{i=1}^{k}\sigma_i(A)u_i v_i^T\in \rt^{m\times n},$ for $k=1,\ldots,\rank (A)$, is the truncated SVD of $A$.


\subsubsection{Review of the Brillinger Transform}\label{bt}

Let $M =  E_{xy}E_{yy}^{- 1}E_{yx}$.
The Brillinger transform (BT)  of $\y$ is represented by  $B = B_1 B_2$ where
$B_1 =  U_{{M},k},$ and $B_2 = U_{{M},k}^T E_{xy} E_{yy}^{-1}$

\subsubsection{Review of Generalized  Brillinger-like Transforms}\label{gbt}

In \cite{tor843}, the generalized Brillinger transforms (GBT1 and GBT2) were developed. The GBT1 is given by ${G_1(\y)=B_1 B_2\y}$, where $ B_1= U_{T_y,k}$, $B_2 = U_{T_y,k}^T E_{x y}E_{yy}^{\dag}$ and $T_y = E_{x y}E_{yy}^{\dag}E_{y x}$. The GBT1 extends the BT to the cases when $E_{yy}$ is singular, i.e., the GBT1 always exists.

The GBT2  generalizes the GBT1 to the form $G_2(\y)$ $=B_1 [B_2\y + B_3\vv]$ where
$B_1= U_{T_u,k},$ $[B_2, B_3]$ $ = U_{T_u,k}^T E_{x u}E_{uu}^{\dag}$
 and $\uu=[\y^T \vv^T]^T\in   L^2(\Omega,\mathbb{R}^{2n})$.  Here, $\vv\in  L^2(\Omega,\mathbb{R}^n)$ is an `auxiliary' signal used to further optimize the transform. If $B_1 B_3\vv=\oo$ then the GBT2 coincides with the GBT1.

\section{The Proposed Approach }\label{new}

\subsection{Generic Transform Structure}\label{nmsmi}

For the  given compression ratio  $c$, we consider  transform  ${\f: \mathbb{R}^{n\times s} \times \mathbb{R}^{q\times s}  \rightarrow \mathbb{R}^{m\times s}}$ represented by
\begin{eqnarray}\label{hq1q2x2}
 \f(Y,V) = D_1 C_1Y + D_2 C_2 Q(Y, V),
\end{eqnarray}
Here,  $D_1,\in\rt^{m\times k_1},$ $ D_2,\in\rt^{m\times k_2},$ $ C_1\in\rt^{k_1\times n},$ $C_2\in\rt^{k_2\times q}$, $k_1 + k_2 = k,$
$V\in  \mathbb{R}^{q\times s}$ is a sample matrix of an `auxiliary' random signal $\vv$ called `injection' and $Q$ is a  transformation of $Y$ and $V$ in matrix $Z\in\mathbb{R}^{q\times s}$, i.e., $Z= Q(Y,V)$.

\subsection{Statement of Problems}\label{nm109}

 Denote $W=[Y^T Z^T]^T$ and write $ \|\cdot \|$ for the Frobenius norm. Find  $D_1,  C_1, D_2, C_2, V$ that solve
 \begin{eqnarray}\label{ei73m}
\min_{V}\min_{{D_1,  C_1}}\min_{{D_2, C_2}} \|X - [D_1C_1 Y + D_2C_2 Q(Y,V)]\|^2,
\end{eqnarray}
and $Q$ which implies
\begin{eqnarray}\label{gg077}
WW^T =\left[ \begin{array}{cc}
                         YY^T  & \oo\\
                         \oo  & ZZ^T
                           \end{array} \right].
\end{eqnarray}
 The condition (\ref{gg077}) leads to the reduction of the numerical load associated with the solution of problem in (\ref{ei73m}). This is because the dimension of $WW^T$ is $(n+q)\times (n+q)$ while the dimensions of  $YY^T$ and  $ZZ^T$ are less, they are $n\times n$ and $q\times q$, respectively. As a result, a solution of problem (\ref{ei73m}) in terms of matrices $YY^T$ and  $ZZ^T$ requires less the numerical load than the solution in terms of matrix  $WW^T$.

The abbreviation MTT will be referred to the multi-term transform based on the solution of the problem in (\ref{ei73m})-(\ref{gg077}).

\begin{remark}
We note that although the source signal $\x$ is not available, the sample matrix $X$ can be represented in terms of matrices $Y$ and sample noise matrix $N$. It is possible, for example, if the observed signal $\y$ is represented in terms of $\x$ and an additive noise (as in Examples  \ref{nm289} and \ref{vnm82} that follow).
\end{remark}

\subsection{Specific features of Problem in (\ref{ei73m})-(\ref{gg077})}\label{nm179}

Interestingly, the problem in (\ref{ei73m})-(\ref{gg077})  can also be treated as a generalization and modification of the blind identification problem considered, in particular, in \cite{Abed1997}.  By the approach in  \cite{Abed1997},  the problem in (\ref{ei73m}) is interpreted differently than that in Section \ref{nm109} above, i.e., $X$ is considered as an available output of the system, $Y$ and $V$ as two inputs, and $D_1,  C_1, D_2, C_2$ make an unknown system function.
Further differences from \cite{Abed1997} are that in (\ref{ei73m}), $Y$ is given, $V$ is unknown and the system function contains more unknowns than those in \cite{Abed1997}. Moreover, matrices  $D_1,  C_1, D_2, C_2$ have special sizes. If we denote ${F_1=D_1C_1}$ and $F_2=D_2C_2$ then $F_1$ and $F_2$ are of ranks $k_1$ and $k_2$, respectively. Therefore, the problem in (\ref{ei73m}) can be reformulated as
  \begin{eqnarray}\label{ei73m1}
\min_{V}\min_{{F_1}}\min_{{F_2}} \|X - [F_1 Y + F_2 Z]\|^2,
\end{eqnarray}
subject to $\rank F_1\leq k_1$, $\rank F_2\leq k_2$ and condition (\ref{gg077}).  By the above reasons, the method in \cite{Abed1997} cannot be applied here. Therefore, in Section \ref{390mn2} below, we propose a new method to solve the problem in (\ref{ei73m})--(\ref{gg077}).

Since $Y$ is available and $V$ is unknown, we call (\ref{ei73m})--(\ref{gg077}) the {\em semi-blind data compression} problem. In terms of \cite{Abed1997}, it can also be called the  {\em semi-blind system identification} problem.

\section{MTT: Solution of problem in (\ref{ei73m}), (\ref{gg077})}\label{390mn2}

Here, we represent an iterative method for finding  optimal matrices $D_1\in \rt^{m\times k_1},$ $  C_1\in\rt^{k_1\times n},$  $D_2\in \rt^{m\times k_2}, C_2\in \rt^{k_2\times q}$, $V$ and transformation $Q$ that solve problem (\ref{ei73m})-(\ref{gg077}).

\subsection{ MTT:  Optimal Semi-Blind Data Compression}\label{bdc1}

First, in (\ref{ei73m}), we determine $Q(Y, V)$  such that
\begin{eqnarray}\label{zzz}
Z = VG \quad\text{where}\quad  G=I- Y^T (Y Y^T)^\dag Y = I- Y^\dag Y.
\end{eqnarray}
Then condition (\ref{gg077}) is true.

As a result, the cost function in (\ref{ei73m}) can be written as
\begin{eqnarray}\label{vbn91}
 &&\hspace*{-12mm}\|X - [D_1C_1 Y + D_2 C_2 Z]\|^2 \nonumber\\
 &&\hspace*{0mm}\hspace*{-2mm} =  \|X - D_1C_1 Y\|^2 +  \|X - D_2 C_2 Z\|^2 -  \|X\|^2.  
\end{eqnarray}

On the basis of representation (\ref{vbn91}), the proposed iterative method consists of the following steps. Let us denote ${S_Y = X Y^T{(YY^T)^{1/2}}^\dag}$ and $T_Y= S_Y S_Y^T =  X Y^T  Y^\dag X^T$.

{\em Step 1:}  For arbitrary $V=V^{(0)}$ and  ${Z^{(0)}} = {V^{(0)}}G$, solve
\begin{eqnarray}\label{ei93m}
\min_{{D_1,  C_1}}\min_{{D_2, C_2}} \|X - [D_1C_1 Y + D_2C_2 Z^{(0)}]\|^2.
\end{eqnarray}

\begin{theorem}\label{nm2}
The solution of problem (\ref{ei93m}) is given by
\begin{eqnarray}\label{30lk2}
&&\hspace*{-12mm}D_1=D^{(0)}_1=U_{T_Y,k_1}, \hspace*{3mm} C_1=C^{(0)}_1=U_{T_Y,k_1}^T  X (Y^T)^\dag,\\
&&\hspace*{-12mm} D_2=D^{(0)}_2\hspace*{-1mm}=U_{T_{Z^{(0)}},k_2},\hspace*{1mm} C_2=C^{(0)}_2\hspace*{-1mm}=\hspace*{-0.5mm}U_{T_{Z^{(0)}},k_2}^T  X ({Z^{(0)}}^T)^\dag. \label{klxm9}
\end{eqnarray}
\end{theorem}

\begin{IEEEproof} Proof is given in Section \ref{app1}.
\end{IEEEproof}

Denote
$$\varepsilon^{(0)} =  \|X - [D^{(0)}_1C^{(0)}_1 Y + D^{(0)}_2C^{(0)}_2 Z^{(0)}]\|^2.$$

{\em Step 2:}
Solve
\begin{equation}\label{xm020}
\min_{V} \|X - [D^{(0)}_1C^{(0)}_1 Y + D^{(0)}_2 C^{(0)}_2 Z]\|^2
\end{equation}
or equivalently, solve
\begin{equation}\label{xm02}
\min_{V} \|X -  D^{(0)}_2C^{(0)}_2 VG\|^2.
\end{equation}
This is because other terms in the RHS of  (\ref{vbn91}) do not depend on $V$.
By \cite{Torokhti2007}, the  minimal norm solution of (\ref{xm02}) is
\begin{equation}\label{xm2m8}
V= V^{(1)} = (D^{(0)}_2C^{(0)}_2)^\dag X G^\dag.
\end{equation}
Denote $Z^{(1)} = V^{(1)}G$ and
\begin{equation}\label{}
\varepsilon^{(1)}_V =  \|X - [D^{(0)}_1C^{(0)}_1 Y + D^{(0)}_2C^{(0)}_2 Z^{(1)}]\|^2.
\end{equation}

{\em Step 3:} Solve
\begin{equation}\label{eitu8m}
\min_{{D_2, C_2}} \|X -  D_2C_2 Z^{(0)}\|^2.
\end{equation}
Similar to (\ref{30lk2}), the solution is given by
\begin{eqnarray}\label{87lk2}
 D_2=D^{(1)}_2=U_{T_Z^{(0)},k_2},\hspace*{1mm} C_2=C^{(1)}_2\hspace*{-1mm}=U_{T_Z^{(0)},k_2}^T  X ([Z^{(0)}]^T)^\dag.
\end{eqnarray}
For the first loop, (\ref{87lk2})  is the same as (\ref{klxm9}) but (\ref{87lk2}) is updated for the next loops (see Step 5 below). Denote
$$
\varepsilon^{(1)}_{D,C} =  \|X - [D^{(0)}_1C^{(0)}_1 Y + D^{(1)}_2C^{(1)}_2 Z^{(0)}]\|^2.
$$

{\em Step 4:}    If $\varepsilon^{(1)}_{V} = \min\{\varepsilon^{(1)}_{D,C}, \varepsilon^{(1)}_{V}\}$ then set $D^{(1)}_j:=D^{(0)}_j$,  $C^{(1)}_j:=C^{(0)}_j$, for $j=1,2$, choose $V^{(1)}$ as in (\ref{xm2m8}) and denote $\varepsilon^{(1)}=\varepsilon^{(1)}_{D,C}$.

If $\varepsilon^{(1)}_{D,C} = \min \{\varepsilon^{(1)}_{D,C}, \varepsilon^{(1)}_{V}\}$ then choose $D^{(1)}_j$ and $C^{(1)}_j$ as above, i.e., $D^{(1)}_j:=D^{(0)}_j$,  $C^{(1)}_j:=C^{(0)}_j$, for $j=1,2$, and set $V^{(1)}:= V^{(0)}$. In this case, denote $\varepsilon^{(1)}=\varepsilon^{(1)}_{V}$.

As a result, the reconstruction of the samples is given by
\begin{eqnarray}\label{xn011}
X^{(1)} = D^{(1)}_1C^{(1)}_1 Y + D^{(1)}_2C^{(1)}_2 Z^{(1)}.
\end{eqnarray}

{\em Step 5:} Repeat steps $2$-$4$ where arbitrary matrix  $V^{(0)}$ is replaced with  $V^{(1)}$. Then in (\ref{eitu8m}), $Z^{(0)}$ is replaced with $Z^{(1)}$.
 The new reconstruction of $X$ is denoted by $X^{(2)}$ which is written similar  to $X^{(1)}$ in (\ref{xn011}).

In general, for $i=1,2,\ldots,$ steps $2$-$4$ are iterated  with updated $V^{(i)}$, $D^{(i)}_2$, $C^{(i)}_2$   and the same $D^{(0)}_1$ and $C^{(0)}_1$. This is  because $D^{(0)}_1$ and $C^{(0)}_1$ do not depend on $V$. The  procedure is continued until a  tolerance $\delta$ is achieved so that $|\varepsilon^{(i+1)}-\varepsilon^{(i)}|\leq \delta$. As a result, compression of $X$ is carried out by $C^{(1)}_1$ and $C^{(i+1)}_2$, and its reconstruction   is carried out by $D^{(1)}_1$ and $D^{(i+1)}_2$.
The associated algorithm is represented below where we denote
\begin{eqnarray}\label{bn01p}
&&\hspace*{-10mm}f(D_1, C_1, D_2, C_2, V)= \|X - [D_1C_1 Y + D_2 C_2 Z]\|^2
\end{eqnarray}
and $f(D_2, C_2,V)=\|X- D_2 C_2VG\|^2.$

{\bf Algorithm MTT:} { Optimal semi-blind data compression by updated optimization of $V^{(i)}$, $D^{(i)}_2$ and $C^{(i)}_2$.}
\smallbreak

Input: $V^{(0)}$ and $\delta\in\rt^+$.

Output:  $D_1^{(0)},$ $C_1^{(0)}$,  $D_2^{(i+1)}$, $C_2^{(i+1)}$ and $V^{(i+1)}$. 

\begin{enumerate}
\item[1.] Solve $\displaystyle\min_{D_1, C_1} \min_{D_2, C_2} f(D_1, C_1, D_2, C_2, V)$;\\
$D_1^{(0)}:= D_1,$  $C_1^{(0)}:=C_1$;  $D_2^{(0)}:= D_2,$  $C_2^{(0)}:=C_2$,\\ $V^{(0)}:=V$;
\item[2.] $\varepsilon^{(0)}=f(D_1^{(0)}, C_1^{(0)}, D^{(0)}_2, C^{(0)}_2, V^{(0)})$;
\item[3.]  for  $i=0,1,2,\ldots$
\item[4.] \hspace{0.5cm}  Solve $\displaystyle\min_{{V}}f(D^{(i)}_2,C^{(i)}_2,{V});$
 \item[5.] \hspace{0.5cm} $\varepsilon^{(i+1)}_{V}=f(D_1^{(0)}, C_1^{(0)}, {D}_2,{C}_2,V);$
\item[6.] \hspace{0.5cm}  Solve $\displaystyle\min_{{D}_2,{C}_2}f({D}_2,{C}_2,V^{(i)});$
    \item[7.] \hspace{0.5cm} $\varepsilon^{(i+1)}_{D, C}=f(D_1^{(0)}, C_1^{(0)}, {D}_2,{C}_2,V^{(i)});$
     \item[8.] \hspace{0.5cm} if $\varepsilon^{(i+1)}_{D, C}\leq \varepsilon^{(i+1)}_{V}$
        \item[9.] \hspace{1cm}   $D^{(i+1)}_2={D}_2$, $C^{(i+1)}_2={C}_2$, $V^{(i+1)}=V^{(i)}$;
        \item[10.] \hspace{1cm}   $\varepsilon^{(i+1)}=\varepsilon^{(i+1)}_{D, C}$;
    \item[11.] \hspace{0.5cm} else
        \item[12.] \hspace{1cm}   $D^{(i+1)}_2=D_2$, $C^{(i+1)}_2=C_2$,   $V^{(i+1)}={V}^{(i)}$;
        \item[13.] \hspace{1cm}   $\varepsilon^{(i+1)}=\varepsilon^{(i+1)}_{V}$;
    \item[14.] \hspace{0.5cm} end
    \item[15.] \hspace{0.5cm} if $|\varepsilon^{(i+1)}-\varepsilon^{(i)}|\leq \delta$
    \item[16.] \hspace{1cm} Stop
    \item[17.] \hspace{0.5cm} end
\item[18.] end
\end{enumerate}

Algorithm MTT is based on the following result.

\begin{theorem}\label{opwm9}
 For $i=0,1,\ldots,$ let $D_1^{(0)},$ $C_1^{(0)},$ $D_2^{(i)},$ $C_2^{(i)}$, $V^{(i)}$ be determined by  Algorithm MTT and let ${\varepsilon^{(i)}=f(D_1^{(0)},C_1^{(0)},D_2^{(i)},C_2^{(i)},V^{(i)})}$ be the associated error. Then the increase in the number of iterations $i$ implies the decrease in the associated error, i.e., $\varepsilon^{(i+1)}\leq \varepsilon^{(i)}.$
\end{theorem}

\begin{IEEEproof} Proof is given in Section \ref{app1}.
\end{IEEEproof}

Convergence of Algorithm MTT  is considered in Section \ref{app1} as well.

\begin{example}\label{ui289}
Let ${\bf y}={\bf s} \circ {\bf x}+10\bxi$ where $\textbf{x}\in L^2(\Omega,\mathbb{R}^{100})$ and $\textbf{s}\in L^2(\Omega,\mathbb{R}^{100})$ are uniformly distributed random signals, and $\xi\in L^2(\Omega,\mathbb{R}^{100})$ is a Gaussian  noise with mean $0$ and variance $1$. Symbol `$\circ$' means the  Hadamard product. Let ${X,S,V\in\mathbb{R}^{100\times 300}}$ be samples of signals  ${\bf x},$ ${\bf s}$  and ${\bf v}$, respectively. Then signal $\y$ is represented by matrix $Y=S\circ X+10N$ where $N \in\mathbb{R}^{100\times 300}$ is a  matrix whose entries are normally distributed with mean $0$ and variance $1$.

For $k_1=k_2=25$ (and $k=k_1+k_2=50$) and a randomly chosen $V=V^{(0)}$, the error associated with   Step 1 of  Algorithm MTT is $\varepsilon^{(0)} = 4.8714$ which is the sane as the error associated with the  GBT2. Let a given tolerance be $\delta = 10^{-5}$. After $19$ iterations of Algorithm MTT, the tolerance is achieved
and the associated error is $\varepsilon^{(19)} = 3.7998$. For  $k=50$, the error associated with  the GBT1 is $\varepsilon(D_1,C_1)=7.5727.$

Thus, the error associated with the MTT is less than those for the GBT1 and GBT2 by $50\%$ and $22\%$, respectively.

In Fig. \ref{ex:fig777}, diagrams of the error associated with the Algorithm MTT are represented, for $100$ experiments. In each experiment, matrix $V=V^{(0)}$ was chosen randomly. Fig. \ref{ex:fig777} illustrate Theorem \ref{opwm9}. That is, the increase in the number of iterations $i$ for updated optimization of $V^{(i)}$, $D^{(i)}_2$, $C^{(i)}_2$ implies the decrease in the associated error.

\begin{figure}[t]
\centering
\begin{tabular}{c}
\includegraphics[scale=0.4]{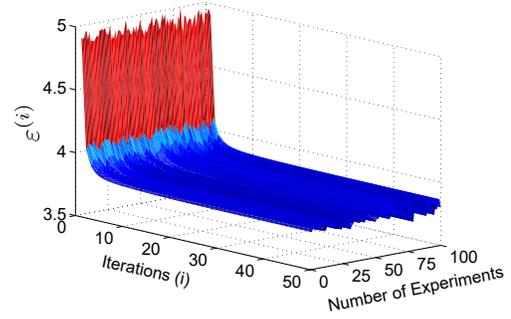}\\
\end{tabular} \caption{Example \ref{ui289}. Diagrams of errors associated with Algorithm MTT.}
  \label{ex:fig777}
\end{figure}

\end{example}

\begin{figure}[h]
    \centering
    \begin{tabular}{cc}
        \includegraphics[scale=0.4]{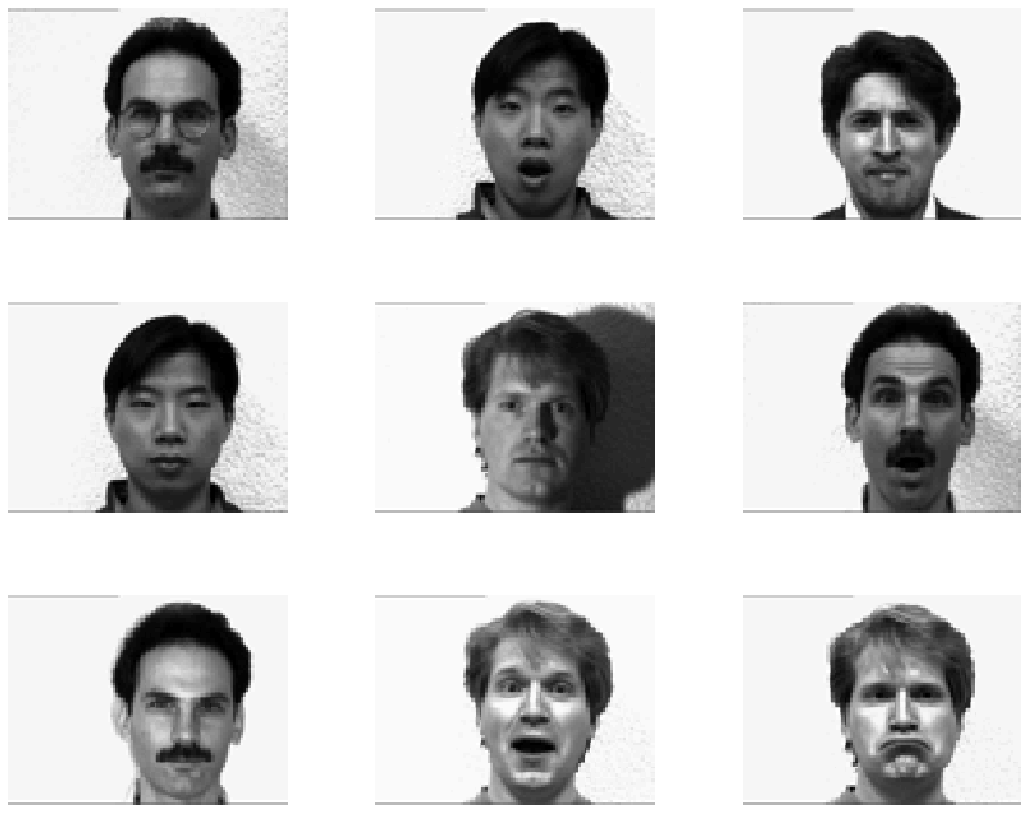}\\  (a)\\  \includegraphics[scale=0.4]{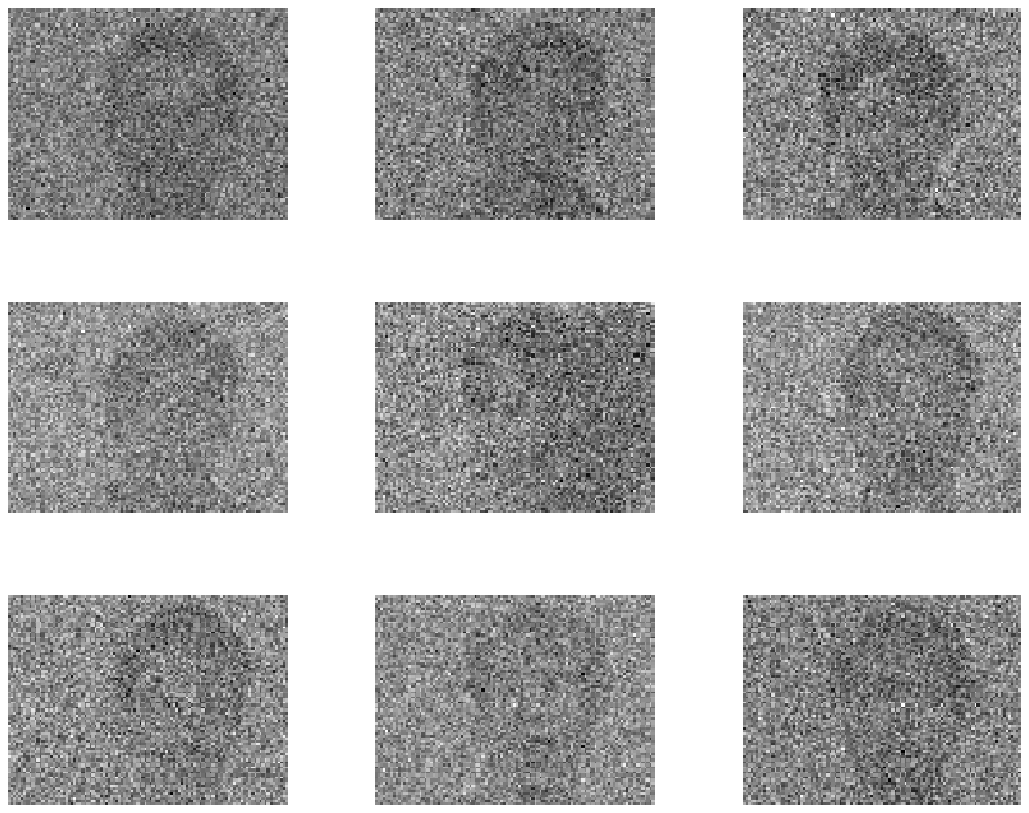}\\
         (b)\\
    \end{tabular}
    \caption{Example \ref{nm289}. (a): Some randomly selected images from the Yale Face Database. (b): Noisy versions of images in (a)}
    \label{example1}
\end{figure}

\begin{example}\label{nm289}
 Here, we illustrate an application of Algorithm MTT to the problem of  compression, filtering and decompression of a set of noisy face images $\mathcal{X}_r=\{\widetilde{X}^{(1)},...,\widetilde{X}^{(r)}\}$ on the basis of a sample of $s\leq r$ images $\mathcal{X}_{t_s}=\{\widetilde{X}^{(t_1)},...,\widetilde{X}^{(t_s)}\}\subset \mathcal{X}_{r}$. The sample $\widetilde{X}^{(t_1)},...,\widetilde{X}^{(t_s)}$ is a permutation of $s$ images from the $r$ images. We choose $r=110$ and $s=55$. Each image  $\widetilde{X}^{(i)}$, for $i=1,\ldots,r$, is simulated by MATLAB as $\widetilde{X}^{(i)}=X^{(i)}+\xi^{(i)}_X$ where $X^{(i)}\in\rt^{81\times 107}$ is a numerical representation of the reference face image taken from the Yale Face Database \cite{Yale_ref}, and  $\xi^{(i)}_X = {\tt randn(81,107)}\in\rt^{81\times 107}$ simulates noise. Some randomly selected face images from the Yale Face Database are given in Fig. \ref{example1} (a). Noisy versions of images in Fig. \ref{example1} (a) are given in Fig. \ref{example1} (b).

Further, let $Y^{(\ell)}\in \mathcal{X}_{r}$, for $\ell=1,...,r$. It is assumed that  $Y^{(\ell)}$ does not necessary belong to the sample $\mathcal{X}_{t_s}$ but is `close' to one of the images in the sample, say $\widetilde{X}^{(\alpha)}$, i.e., $\widetilde{X}^{(\alpha)}$, for $\alpha= t_1,\ldots,t_s$, is such that, for all element in $\mathcal{X}_{t_s}$,
$$
\widetilde{X}^{(\alpha)} \in \text{arg}\min_{\nu=t_1,\ldots,t_s} \|\widetilde{X}^{(\nu)} - Y^{(\ell)}\|^2\leq \delta,
$$
for a given $\delta \geq 0.$ We wish to formulate the problem under consideration in terms of the problem in (\ref{ei73m}). To this end, we denote
$$
X=[X^{(t_1)},...,X^{(t_s)}]\qa
Y=[\widetilde{X}^{(t_1)},...,\widetilde{X}^{(t_s)}]
$$
where $X,Y\in \mathbb{R}^{81\times (107\times s)}$.  Matrix $V$ (see (\ref{ei73m})) is simulated as
$V=[V^{(1)},...,V^{(g)}]\in \mathbb{R}^{81\times q}$ where $q=107\times g$ and entries of each $V^{(i)}$, for $i=1,\ldots,g$, are uniformly distributed. In this notation, the problem under consideration is formulated as in (\ref{ei73m}), i.e., we wish to find $D_1,  C_1, D_2, C_2$ and $V$ that solve
  \begin{eqnarray}\label{jk73m}
\min_{V}\min_{{D_1,  C_1}}\min_{{D_2, C_2}} \|X - [D_1C_1 Y + D_2C_2 Z]\|^2.
\end{eqnarray}
Matrix $Z$ is represented it in the block form as $Z=[Z^{(t_1)},...,Z^{(t_s)}]$ where $Z^{(j)}\in\rt^{k_2\times 107}$, for $j=1,\ldots, s$.

Matrices $D_1,  C_1, D_2, C_2$ and $V$ are determined by Algorithm MTT. Then an estimate of  $X^{(\alpha)}$, for $\alpha=1,\ldots, r$,  is given by $\widehat{X}^{(\alpha)} = D_1C_1 Y^{(\ell)} + D_2C_2 Z^{(\ell)}$.  The estimate is represented, for $k_1=20$, $k_2=20$ and $k=k_1+k_2=40$, in  Fig. \ref{example2} (e), for $V$ with arbitrarily uniformly distributed entries,  and in Fig. \ref{example2} (f) with optimal $V$ using $10$ iterations.

The GBT1 and GBT2 \cite{tor843} have also been applied to the problem in consideration, for the same $k=40$. The mean square errors (MSEs) are represented in Fig. \ref{example2}  and Table 1 where the MSE means the mean square error. The obtained numerical results clearly demonstrate the advantages of the Algorithm MTT. The MSE and SSIM associated with the MTT are much better than those for the other methods.

\begin{figure}[t]
    \centering
    \begin{tabular}{ccc}
          \includegraphics[scale=0.4]{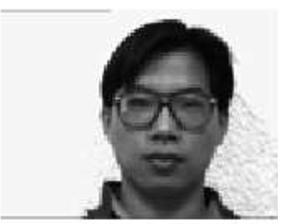} &  \includegraphics[scale=0.4]{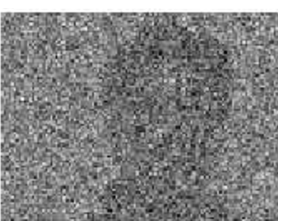} &  \includegraphics[scale=0.4]{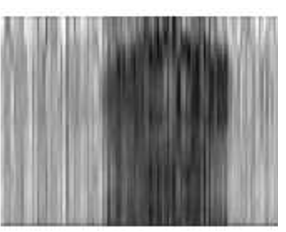} \\
         (a) & (b) & (c)\\
    \end{tabular}
    \begin{tabular}{cc}
          \includegraphics[scale=0.4]{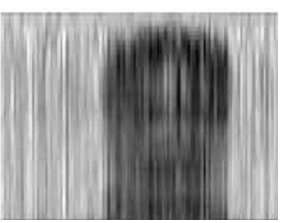} &   \includegraphics[scale=0.4]{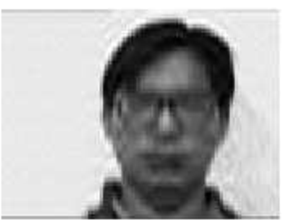} \\
         (d) & (e) \\
    \end{tabular}
    \caption{  Example \ref{nm289}:  Reference image ${X}^{(\alpha)}$ (a);   Observed image $\widetilde{Y}^{(\ell)}$ (b); Estimates of ${X}^{(\alpha)}$ by the GBT1 (c);  GBT2  (d);  MTT (e).}
    \label{example2}
\end{figure}

    \centering
    \begin{tabular}{|c|c|}
     \multicolumn{2}{c}{Table 1}\\
      \hline
     Transforms  & MSEs \\
      \hline
      GBT1 & $3.50\times 10^2$ \\
      \hline
      GBT2 & $2.80\times 10^2$ \\
      \hline
      MTT & $1.80\times 10^{1}$ \\
      \hline
    \end{tabular}

\end{example}

\subsection{Particular Feature of MTT}\label{nma9}

On the basis of results presented in \cite{Aa2017}, it can be shown that the MTT converges faster if the dimension $q$ of matrix $V$ increases. Here, we illustrate this feature as follows.

\begin{example}\label{vnm82}
 Let ${\bf y}={\bf x}+\bxi$ where $\x\in L^2(\Omega,\mathbb{R}^{50})$ and $\bxi\in L^2(\Omega,\mathbb{R}^{50})$  are uniformly distributed source signal and  white noise, respectively. Let ${\bf v}\in L^2(\Omega,\mathbb{R}^{q})$ be a  Gaussian random signal.
 Noise $\bxi$ is  uncorrelated with ${\bf x}$. Vectors $\y$, $\bxi$ and $\vv$ are represented by samples $Y\in\rt^{50\times 100}$, $\Xi\in\rt^{50\times 100}$ and $V\in\rt^{q\times 100}$. Then
$
E_{yy}=\frac{1}{100}YY^T, \quad E_{xx}=E_{yy}-E_{\xi\xi} \qa E_{xy}=E_{xx}
$
where
$E_{\xi\xi}=\sigma^2 I$ and $\sigma=2$.

In Fig. \ref{ex:fig107}, for $k_1=k_2=12$, diagrams of the error associated with Step 1 of the MTT  versus dimension $q$ of matrix $V$ are given, for $100$ experiments. The error decreases when $q$ increases. In each experiment, matrix $Y$ was chosen randomly.

\begin{figure}[h!]
\centering
\includegraphics[scale=0.4]{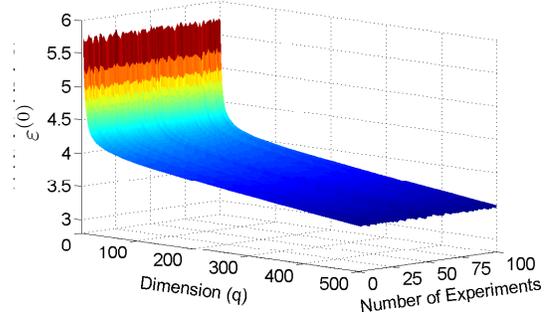}
  \caption{Example \ref{vnm82}. Diagrams of  the error associated with Step 1 of the MTT versus dimension $q$ of matrix $V$.}
  \label{ex:fig107}
\end{figure}

\end{example}

\vspace*{-3mm}
\section{Appendix}\label{app1}

\subsubsection{Proof of Theorem \ref{nm2}} Denote $W:=W^{(0)} = [Y^T {Z^{(0)}}^T]^T$.
\begin{eqnarray}\label{kbnb4}
&&\hspace*{-15mm}\mbox{Then}\hspace*{3mm} \|X - [D_1C_1 Y + D_2C_2 Z^{(0)}]\|^2 \nonumber\\
&&\hspace*{-8mm}=\|(XX^T)^{1/2}\|^2 - \|XW^T{(WW^T)^{1/2}}^\dag\|^2 \\
&& + \|XW^T(WW^T)^{1/2} - [F_1 F_2](WW^T)^{1/2}\|^2.\nonumber
\end{eqnarray}
Further,
$
{(WW^T)^{1/2}}^\dag=\left[ \begin{array}{cc}
                         {(YY^T)^{1/2}}^\dag  & \oo\\
                         \oo  & {(ZZ^T)^{1/2}}^\dag
                           \end{array} \right].
$ Then
\begin{eqnarray*}\
&&\hspace*{-13mm}\|XW^T(WW^T)^{1/2} - [F_1 F_2](WW^T)^{1/2}\|^2 \\
&& = \|S_Y - F_1{(YY^T)^{1/2}}\|^2 + \|S_Z - F_2{(ZZ^T)^{1/2}}\|^2.
\end{eqnarray*}
Further, let us denote by $\rt(m, n, k)$ the set of all $m \times n$ matrices of rank  at most $k$.
The minimal norm solutions to the problems
$
\displaystyle \min_{F_1\in \mathbb R(m, n, k_1) } \|S_Y - F_1{(YY^T)^{1/2}} \|^2
$
and
$
\displaystyle\min_{F_2\in \mathbb R(m, n, k_2) } \|S_Z - F_2{(ZZ^T)^{1/2}}\|^2
$
are given  \cite{Torokhti2007} by
\begin{eqnarray}\label{sn021m}
F_1 = \left[ S_Y\right]_{k_1} {(YY^T)^{1/2}}^\dag \qa
 F_2 = \left[ S_Z\right]_{k_2}  {(ZZ^T)^{1/2}}^\dag,
\end{eqnarray}
respectively. Then (\ref{30lk2}), (\ref{klxm9}) follow from (\ref{sn021m})  on the basis of Fact 1 in \cite{tor843}.


\subsubsection{Proof of Theorem \ref{opwm9}} To prove Theorem \ref{opwm9} we need some preliminaries as follows.

\begin{definition}
Let $\mathcal M$ be a metric space and let $\mathbb{B}(\alpha,\epsilon) \subset \mathcal M$ be an open ball with center $\alpha$ and radius $\epsilon$.
 Point $\alpha$  is called  a {\em cluster point} \cite{Amann2006} of a sequence $\{x_n\} \subset \mathcal M$ iff for each $\epsilon>0$ and $m\in\mathbb{N}$, there is some $n\geq m$ such that $x_n\in\mathbb{B}(\alpha,\epsilon)$, for all $n\geq m$.
\end{definition}
We denote
$
F_j = D_j C_j,$ $F^{(i)}_j = D^{(i)}_jC^{(i)}_j,
$
 $F=\{F_1, F_2\}$ and $F^{(i)}=\{F^{(i)}_1, F^{(i)}_2\}$ where $j=1,2$ and $i=0,1,\ldots$. Then
 \begin{eqnarray*}
&&\hspace*{-11mm}f(F, V) = f(F_1,F_2,V)\\
&&\hspace*{1.5mm}=f(D_1, C_1, D_2, C_2, V)=\|X - F_1 Y - F_2VG\|^2
 \end{eqnarray*}
 where $f(D_1, C_1, D_2, C_2, V)$ is defined by (\ref{bn01p}).

Let us now define compact sets $K_1$ and $K_2$  such that
\begin{eqnarray}\label{k1234}
K_1=\{F:0\leq f(F,V^{(0)})\leq f(F^{(0)},V^{(0)})\}
\end{eqnarray}
\begin{eqnarray}\label{k1235}
K_2=\{V:0\leq f(F^{(0)},V)\leq f(F^{(0)},V^{(0)})\}.
\end{eqnarray}

\begin{definition}
Let $S^{*}=(F^{*},V^{*})$ be a cluster point of sequence $S^{(j)}=(F^{(j)},V^{(j)})$, for $j=1,2,\ldots.$ Point $\overline{F}_{V^*}$ is called a {\em best response} to $F$ if
$
\displaystyle\overline{F}_{V^*}\in\arg\min_{F\in K_1}f(F,V^*).
$
\end{definition}
\medbreak

The proof of Theorem \ref{opwm9}  is as follows.

\begin{IEEEproof}
For given $F_1^{(0)}=D_1^{(0)}C_1^{(0)}$, $F_2^{(i)}=D_2^{(i)}C_2^{(i)}$ and $V^{(i)}$, we determine $\varepsilon^{(i)}=f({F}^{(0)}_1,F_2^{(i)},V^{(i)}).$ Then $\varepsilon^{(i+1)}$ is determined with the following steps.
By  Algorithm MTT, the best responses are computed as follows:

 Given $F^{(0)}_1$ and $F^{(i)}_2$, compute $\widetilde{V}$ which is the best response to $V$ for $f(F^{(0)}_1,F^{(i)}_2,V)$.
 Given $F^{(0)}_1$ and $V^{(i)}$, compute $\widetilde{F}_2$ of $\rank \leq k_2$ which is the best response to $F_2$ for $f(F^{(0)}_1,F_2,V^{(i)})$.

Let us now define
$
\varepsilon^{(i)}_{F_2}=f(F^{(0)}_1,\widetilde{F}_2,V^{(i)})$ and $\varepsilon^{(i)}_{V}=f(F^{(0)}_1,F^{(i)}_2,\widetilde{V}).
 $
 If $\varepsilon^{(i)}_{F_2}=\min\{\varepsilon^{(i)}_{F_2}, \varepsilon^{(i)}_{V}\},$ then $$\varepsilon^{(i+1)}=\varepsilon^{(i)}_{F_2}=f(F^{(0)}_1,\widetilde{F}_2,V^{(i)}).$$
  At the same time, for given $F^{(0)}_1$ and $V^{(i)}$, the best response to $F_2$, for $f(F^{(0)}_1,F_2,V^{(i)})$, is $\widetilde{F}_2$. Then, for all $F_2$,
  $$
  f(F^{(0)}_1,\widetilde{F}_2,V^{(i)})\leq f(F^{(0)}_1,F_2,V^{(i)}).
  $$
In particular, $f(F^{(0)}_1,\widetilde{F}_2,V^{(i)})\leq f(F^{(0)}_1,F_2^{(i)},V^{(i)}).$ Then
$ 
    \varepsilon^{(i+1)}\leq \varepsilon^{(i)}.
$ 
If $\varepsilon^{(i)}_{V}=\min\{\varepsilon^{(i)}_{F_2}, \varepsilon^{(i)}_{V}\},$ then
$$
\varepsilon^{(i+1)}=\varepsilon^{(i)}_{V}=f(F^{(0)}_1,F^{(i)}_2,\widetilde{V}).
 $$
 Further, given $F^{(0)}_1$ and $F^{(i)}_2$, the best response to $V$, for $f(F^{(0)}_1,F_2^{(i)},V),$ is $\widetilde{V}$.
    Then, for all $V$, we obtain
    $$
    f(F^{(0)}_1,F^{(i)}_2,\widetilde{V})\leq f(F^{(0)}_1,F_2^{(i)},V).
    $$
    In particular,
    $$
    f(F^{(0)}_1F^{(i)}_2,\widetilde{V})\leq f(F^{(0)}_1,F_2^{(i)},V^{(i)}).
    $$
    i.e.,
  $ 
   \varepsilon^{(i+1)}\leq \varepsilon^{(i)}.
  $ 
Then the statement of Theorem \ref{opwm9} follows.
 \end{IEEEproof}

\subsubsection{Convergence of Algorithm MTT}
\begin{theorem}\label{opw9}
Let $K_1$ and $K_2$ be compact sets defined by (\ref{k1234}) and (\ref{k1235}), respectively, and
$F\in K_1$ and $V\in K_2$. Let $F^{(j)}$ and $V^{(j)}$ be determined by Algorithm MTT, for $j=0,1,\ldots.$
Then any cluster point of  sequence $S^{(j)}=(F^{(j)},V^{(j)})$, say $S^{*}=(F^{*},V^{*})$, is a coordinate-wise minimum point of
$ f(F,V),$
i.e.,
$$
F^{*}\in\arg\min_{F\in K_1} f(F,V^*),\;\;\;\;V^{*}\in\arg\min_{V\in K_2} f(F^*,V).
$$
\end{theorem}

 \begin{IEEEproof} Since each $K_j$ is compact, then there is a subsequence $S^{(j_t)}=(F^{(j_t)},V^{(j_t)})$ such that $S^{(j_t)}\rightarrow S^{*}$ when $t\rightarrow \infty$.
Consider  entry $F^{(j_t)}$ of $S^{(j_t)}$. Let $\overline{F}_{V^*}$ and  $\overline{F}_{V^{(j_t)}}$ be best responses to $F$ associated with $\overline{F}_{V^*}$ and $V^{(j_t)}$, respectively.
Then we have
\begin{eqnarray*}
&&\hspace*{-10mm}f(\overline{F}_{V^*},V^{(j_t)})\geq   f(\overline{F}_{V^{(j_t)}},V^{(j_t)})\geq  f(F^{(j_t+1)},V^{(j_t+1)})\\
                                 &&\hspace*{20mm}\geq  f(F^{(j_{t+1})},V^{(j_{t+1})})
\end{eqnarray*}
By continuity, $f(\overline{F}_{V^*},V^{*})\geq f(F^{*},V^{*})$ as $t\rightarrow \infty$.
 It implies that latter should hold as an equality, since the inequality is true by the definition of the best response $\overline{F}_{V^*}$. Thus, $F^{*}$ is the best response for $V^*$, or equivalently, $\overline{F}^*$ is the solution for the problem
$$\displaystyle\arg\min_{F\in K_1} f(F,V^*).$$

The proof  is similar if we consider entry $V^{(j_t)}$ of $S^{(j_t)}$.
 \end{IEEEproof}

\bibliographystyle{IEEEtran}
\bibliography{Aaachen_2_2017}

\end{document}